\documentclass[12pt]{amsart}

\usepackage{amsfonts,amsmath,amssymb,amsthm}

\newcommand{\Cc}{\mathbb{C}} 
\newcommand{\Pp}{\mathbb{P}}

\newcommand{\Nn}{\mathbb{N}}
\newcommand{\Zz}{\mathbb{Z}}
\newcommand{\Qq}{\mathbb{Q}}

\newcommand{\BB}{\mathcal{B}} 
\newcommand{\CC}{\mathcal{C}} 
\newcommand {\Sing}{\mathop{\mathrm{Sing}}\nolimits}
\newcommand {\Aut}{\mathop{\mathrm{Aut}}\nolimits}
\newcommand{\defi}[1]{\emph{#1}}

\renewcommand{\epsilon}{\varepsilon}
\renewcommand{\le}{\leqslant}
\renewcommand{\ge}{\geqslant}

{\theoremstyle{plain}
\newtheorem{theorem}{Theorem}    % theorem with number
\newtheorem*{theorem*}{Theorem}    % theorem with number
\newtheorem{lemma}[theorem]{Lemma}       % lemma with number
}
{\theoremstyle{remark}
\newtheorem*{remark*}{Remark}  % theorem without number
\newtheorem{example}[theorem]{Example}
}

\begin{document}
%%%%%%%%%%%%%%%%%%%%%%%%%%%%%%%%%%%%%%%%%%%%%%%%%%%%%%%%%%%%%%%%%
%%%%%%%%%%%%%%%%%%%%%%%%%%%%%%%%%%%%%%%%%%%%%%%%%%%%%%%%%%%%%%%%%

\title{Integral points on generic fibers}

\author{Arnaud Bodin}
\address{Laboratoire Paul Painlev\'e, Math\'ematiques,
Universit\'e de Lille 1,  59655 Villeneuve d'Ascq, France.}
\email{Arnaud.Bodin@math.univ-lille1.fr}

\date{\today}

\begin{abstract}
Let $P(x,y)$ be a rational polynomial and $k\in \Qq$ be a generic value.
If the curve $(P(x,y)=k)$ is irreducible and admits an infinite number of points whose coordinates are integers
then there exist algebraic automorphisms  that send $P(x,y)$ to the polynomial $x$
or to $x^2-dy^2$, $d\in \Nn$. 
Moreover for such curves (and others) we give a sharp bound for the number
of integral points $(x,y)$ with $x$ and $y$ bounded.
\end{abstract}

% MSC2000 14H50, 14R05

\maketitle

%%%%%%%%%%%%%%%%%%%%%%%%%%%%%%%%%%%%%%%%%%%%%%%%%%%%%%%%%%%%%%%%%
\section{Introduction}

Let $P \in \Qq[x,y]$ be a polynomial and $\CC = (P(x,y)=0) \subset \Cc^2$ be the corresponding algebraic curve.
On old and famous result is the following:
\begin{theorem*}[Siegel's theorem]
Suppose that $\CC$ is irreducible.
If the number of integral points $\CC \cap \Zz^2$ is infinite then
$\CC$ is a rational curve.  
\end{theorem*}

Our first goal is to prove a stronger version of Siegel's theorem
for curve defined by an equation $\CC = (P(x,y)=k)$ where $k$ is a generic value.
More precisely there exists a finite set $\BB$
such that the topology of the complex plane curve
$(P(x,y)=k) \subset \Cc^2$ is independent of $k\notin \BB$.
We say that $k\notin \BB$ is a \defi{generic value}.

\begin{theorem}
\label{th:coord}
Let $P \in \Qq[x,y]$ and let $k \in \Qq\setminus \BB$ be a generic value.
Suppose that the algebraic curve $\CC  = (P(x,y)=k)$ is irreducible.
If $\CC$ contains an infinite number of integral
points $(m,n) \in \Zz^2$ then there exists an algebraic automorphism $\Phi \in \Aut \Qq^2$ such 
$$P \circ \Phi (x,y) = x \quad \text{ or } \quad  P \circ \Phi (x,y) = \alpha(x^2-dy^2)+\beta,$$
where $d\in \Nn^*$ is a non-square and $\alpha\in \Qq^*$, $\beta \in \Qq$.
\end{theorem}

In particular the curve $\CC  = (P(x,y)=k)$ is diffeomorphic to a line $(x=0)$, in which case
the set $\BB$ is empty or to an hyperbola $x^2-dy^2=1$ in which case $\BB$ is a singleton.

Theorem~\ref{th:coord} can be seen as an arithmetic version of the Abhyankar-Moh-Suzuki theorem
\cite{AM} and in fact we use this result.
It can also be seen as a strong version of a result of Nguyen Van Chau \cite{Nc}
concerning counter-examples to the Jacobian conjecture.

For example, let $P(x,y) = x-y^d$. The curve $\CC = (x - y^d=0)$ has infinitely 
many integral points of type $(n^d,n)$, $n\in \Zz$.
And for the algebraic automorphism $\Phi(x,y) = (x+y^d,y)$
we have $P \circ \Phi(x,y) = x$.
In particular the curve $\CC$ is send (by $\Phi^{-1}$) to
$(x=0)$. 
As pointing out by Kevin Buzzard the second case corresponds to Pell's equation 
and for example $(x^2-2y^2=1)$ admits an infinite number of integral points.
Of course a kind of reciprocal of Theorem~\ref{th:coord} is true.
Let $Q_1(x,y) = x$, (resp.\ $Q_2(x,y) = x^2-dy^2$) and $\Phi\in \Aut \Qq^2$ 
whose inverse $\Phi^{-1}$ as integral coefficients.
If we set $P_1 = Q_1 \circ \Phi$ (resp.\ $P_2 = Q_2 \circ \Phi$) then
the curve $(P_1=0)$  (resp.\ $(P_2=1)$) has infinitely many integral points.

For non-generic values the result is not true, for example let 
$P(x,y) = x^2 - y^3$ and $\CC = (x^2-y^3 = 0)$.
The integral points $(n^3,n^2)$, $n\in \Zz$ belongs to $\CC$,
but as $\CC$ is singular it cannot be algebraically equivalent to a line. 

\bigskip
\bigskip

We will apply Theorem \ref{th:coord} to obtain new bounds for the number of 
integral points on algebraic curves.
Let $\CC = (P(x,y)=0)$ be an algebraic curve, and let $d = \deg P$.
Let 
$$N(\CC, B) = \# \left\lbrace (x,y) \in \CC \cap \Zz^2 \mid |x| \le B \text{ and } 
|y| \le B\right\rbrace.$$
For all irreducible curves Heath-Brown \cite{Hb} proved that $N(\CC,B) \le C_d B^{\frac 1 d + \epsilon}$
for some constant $C_d$. By making the proof of \cite{Hb} explicit
Walkowiak obtains in \cite{Wa}:
\begin{theorem*}[Explicit Heath-Brown's theorem]
For all irreducible curve $\CC$ of degree $d$ and all $B>0$: 
 $$N(\CC,B) \le 2^{48}d^8 \ln(B)^5 B^{\frac 1 d}.$$
\end{theorem*}

The term $B^{\frac 1 d}$ in Theorem \ref{th:numb} is sharp but the term $2^{48}d^8 \ln(B)^5$ 
is far from being optimal.
For curves $\CC$ as in Theorem \ref{th:coord} we will give sharp bounds
for $N(\CC,B)$. 
First of all if $\CC = (P=k)$ and the polynomial $P$ is algebraically equivalent to $x^2-dy^2$
it is know \cite[p.~135]{Se} that there exists $C >0$ and $n \ge 1$ such that
$$N(\CC,B) \le C \cdot \ln(B)^n.$$
Of course it implies $N(\CC,B) \le B^{\frac 1 d}$ for sufficiently large $B$
and we shall omit this case. 

So if $P$ is not algebraically equivalent to $x^2-dy^2$ then
 as a corollary of Theorem \ref{th:coord}, 
such a curve $\CC$ admits a parametrisation by polynomials:
let $(p(t),q(t))$ be a parametrisation of $\CC$
with rational coefficients. Moreover $\deg P$ is equal to $\deg p$ or $\deg q$.
We suppose $\deg P = \deg p$ and write 
$$ p(t) = \frac{1}{b}(a_d t^d + a_{d-1}t^{d-1}+\cdots+a_0).$$
where $a_0,\ldots, a_d,b \in \Zz$ and $\gcd(a_0,\ldots, a_d,b) = 1$.

\begin{theorem}
\label{th:numb}
Let $\CC$ be as in Theorem \ref{th:coord}. 
Then the number $N(\CC, B)$  of integral points on $\CC$ bounded by $B$, verifies:
$$ \forall \epsilon > 0 \quad \exists B_0 \ge 0  \quad \forall B \ge B_0 \qquad 
N(\CC,B) \le 2a_d^{1-\frac 1 d} b^{\frac 1 d} B^{\frac 1 d} + 1 + \epsilon.$$
\end{theorem}

In fact by adding a term $\frac{(d-1)(d-2)}{2}$ we prove such a bound 
for curves $\CC$ parametrised by polynomials $(p(t),q(t))$.

For example if  $a_d = 1$ and $b=1$ in Theorem \ref{th:numb}
we get for a sufficiently large $B$:
$$N(\CC,B) \le 2B^{\frac 1 d} + 1 + \epsilon.$$
It implies for $\epsilon = \frac 12$ and all sufficiently large $B$ such that $B^{\frac{1}{d}}$
is an integer that we get 
$$N(\CC, B) \le 2B^{\frac 1 d} + 1.$$
For instance we have a parametrisation of $\CC = (x-y^d=0)$
by $(t^d,t)$. If $B^{\frac{1}{d}}$ is an integer we get 
$N(\CC,B) = 2B^{\frac 1 d} + 1$.
Moreover the ``$\epsilon$'' term is necessary as
the example $\CC = (x-y^d=1)$ proves, see Example \ref{ex:epsilon}.

\bigskip
\bigskip

\emph{Acknowledgements:} This work has been done during a visit at the University
of Zaragoza, the author wishes to thank people of the  Department of Mathematics and 
especially Enrique Artal Bartolo for hospitality. I also thank Pierre D\`ebes
for his support.

%%%%%%%%%%%%%%%%%%%%%%%%%%%%%%%%%%%%%%%%%%%%%%%%%%%%%%%%%%%%%%%%%
\section{Parametrisation}
\label{sec:param}

\subsection{Topology of polynomials}

By a result of Thom for a polynomial $P\in \Qq[x,y]$ seen as a map $P : \Cc^2 \rightarrow \Cc$
there exists a finite set $\BB \subset \Cc$ such that 
$$P : P^{-1}(\Cc \setminus \BB) \longrightarrow \Cc \setminus \BB$$
is a topological locally trivial fibration.
A value $k \notin \BB$ is a \defi{generic value}.

For example we have the following characterization of the generic values:
the Euler characteristic of the complex plane curve 
$(P(x,y)=k)\subset \Cc^2$ is independent of $k\notin \BB$ and jumps if and only if $k\in \BB$.

Of course the image by $P$ of a singular point is not a generic value, but for example 
$P(x,y) = x(xy-1)$ has no singular points while $\BB = \{ 0 \}$.
Then if $k \notin \BB$ is a generic value the plane algebraic curve
$\CC = (P(x,y)=k) \subset \Cc^2$ does not have singular points.
Moreover if $\CC$ is connected then $\CC$ is irreducible.

The connectedness of a generic fiber $\CC$ is equivalent of $P(x,y)$ being 
non-composite \cite{Ar}.  We recall that $P(x,y)$ is \defi{composite} 
if there exist $h \in \Cc[t]$, $\deg h \ge 2$, and $Q \in \Cc[x,y]$ 
such that $P(x,y) = h \circ Q(x,y)$. By \cite[Theorem 7]{Ay} we even 
can choose $h$ and $Q$ with rational coefficients.
Consequently it has been noticed by Janusz Gwozdziewicz that the hypothesis ``$\CC$ is irreducible''
in Theorem \ref{th:coord} can be removed. In that case the conclusion becomes
$P \circ \Phi (x,y) = h(x)$ or $P \circ \Phi (x,y) = h(x^2-dy^2)$,
where $h \in \Qq[t]$ is a one-variable polynomial of positive degree.

\subsection{Algebraic automorphisms}

For $K = \Qq$ or $K=\Cc$ an \defi{algebraic automorphism} $\Phi \in \Aut K^2$
is a polynomial map $\Phi : K^2 \longrightarrow K^2$, invertible, such that the inverse is also a
polynomial map. The polynomials $P,Q \in K[x,y]$ are \defi{algebraically equivalent} if there exists
$\Phi \in \Aut K^2$ such that $Q = P \circ \Phi$. And in fact 
such $P$ and $Q$ have the same topological and algebraic properties.

\subsection{Siegel's theorem}

By Siegel's theorem an irreducible plane algebraic curve $\CC$ with an infinite 
number of integral points $(x,y) \in \CC \cap \Zz^2$ is a rational curve (i.e.\ the genus is zero).

\subsection{Parametrisation}

As $\CC = (P(x,y)=k)$ is rational it admits a parametrisation by rational fractions.
In order to deals with points at infinity and special parameters we compactify the situation.

Let $\bar P(x,y,z)$ be the homogenisation of $P(x,y)$ with $d = \deg \bar P = \deg P$.
Then $\bar \CC = (P(x,y,z)-k z^d=0) \subset \Pp^2$ is the closure of $\CC$.

As $\bar \CC$ is rational there exists a birational map 
$\phi : \Pp^1 \longrightarrow \bar \CC$ defined by
$$\phi (t,s) = (\bar p(t,s),\bar q(t,s), \bar r(t,s))$$
where $\bar p, \bar q, \bar r$ are homogeneous polynomials of the same degree in $\Zz[t,s]$,
without common non-constant factor.

We will need some facts about parametrisations (see \cite{PV}).
\begin{lemma}
\label{lem:param}
For any such parametrisation:
\begin{enumerate}
  \item \label{it:l1} $\phi$ is a morphism (it is well-defined everywhere);
  \item \label{it:l2} $\phi$ is surjective;
  \item \label{it:l3} $\deg \bar p = \deg \bar q = \deg \bar r = \deg P = d$;
  \item \label{it:l4} The birational inverse $\psi$ of $\phi$ is well-defined 
away from the singular points of $\bar \CC$;
  \item \label{it:l5} If $(x,y) \in \CC \cap \Zz^2$ is a non-singular point then
there exists $(t,s) \in \Zz^2$ such that $\phi(t,s) = (x:y:1)$;
\end{enumerate}
\end{lemma}

\begin{proof}
We just sketch the proofs and we refer to \cite[Lemma 2.1]{PV} for details and references:
the fact that $\bar p, \bar q, \bar r$ have no common factor implies that 
$\phi$ is well-defined. Then $\phi(\Pp^1)$ is dense in $\bar \CC$
and hence is equal to $\bar \CC$. 
The birational inverse $\psi$ of $\phi$ is well-defined excepted at the singular points of $\bar \CC$,
see Fulton \cite[p.160]{Fu}. Then a non-singular integral point of $\bar \CC$ is send by $\psi$ to
a point with rational coordinates.
We will not need item (\ref{it:l3}) and we refer to \cite{PV}.
\end{proof}

\subsection{Maillet's theorem}

An old result of Maillet \cite{Ma} gives strong restrictions for the parametrisations.
\begin{theorem}
\label{th:maillet}
If $\CC = (P(x,y)=k) \subset \Cc^2$ has an infinite number of integral points then there
exists a parametrisation $\phi$ of $\bar \CC$ given by
 $$\phi (t,s) = (\bar p(t,s),\bar q(t,s), \bar r(t,s))$$
as before with
$$ \bar r(t,s) = at^d \quad \text{ or } \quad \bar r(s,t) = a(\alpha t^2 + \beta ts + \gamma s^2)^{d/2}$$
where $a,\alpha,\beta,\gamma \in \Zz$ and $\beta^2-4\alpha\gamma >0$.
\end{theorem}

\subsection{Topology of $\CC$}

\begin{lemma}
\label{lem:topo}
If $k$ is a generic value and $\CC = (P(x,y)=k)$ has an infinite number of integral points then 
$\CC$ is homeomorphic to $\Cc$ or to $\Cc^*$.
\end{lemma}

In fact the homeomorphisms are 
diffeomorphisms and we only need $k$ to be a non-singular value of $P(x,y)$. In Lemma \ref{lem:neumann} we will exclude the case $\Cc^*$.

\begin{proof}
Let $L_\infty = (z=0)$ of $\Pp^2$ be the line at infinity.
What are the parameters
$(t:s)$ that give the points at infinity $\bar \CC \cap L_\infty$? 
The points at infinity correspond to the parameters $(t:s)$ such that
$\bar r(s,t)=0$. 

So if $\phi$ is the parametrisation given by Maillet's theorem
then for $\bar r(t,s) = at^d$ then $(t:s)= (0:1)$ and there is one point $P \in \bar \CC \cap L_\infty$.
For $\bar r(s,t) = a(\alpha t^2 + \beta ts + \gamma s^2)^{d/2}$ 
then $(t:s) = (\tau_1:1)$ or $(t,s) = (\tau_2:1)$
where $\tau_1,\tau_2$ are the roots of $\alpha t^2 + \beta t + \gamma$ (say $\alpha \neq 0$);
there are two points $P_1,P_2 \in \bar \CC \cap L_\infty$.

From Lemma \ref{lem:param} we know that
$\phi$ is morphism that induces a bijective map
onto $\bar \CC \setminus \Sing \bar \CC$.
The only singular points of $\bar \CC$ are on the line at infinity $L_\infty$ 
because $\CC$ is a non-singular affine curve.
Then in the case $\bar r(s,t) = at^d$, the restriction
$\phi_| : \Pp^1 \setminus \{ (0:1) \} \longrightarrow \bar \CC \setminus \{ P \}$
is a bijective map. Moreover it is an homeomorphism. But $\Pp^1 \setminus \{ (0:1) \} = \Cc$ 
and $\bar \CC \setminus \{ P \} = \CC$. Then $\CC$ is homeomorphic to $\Cc$.

In the case $\bar r(s,t) = a(\alpha t^2 + \beta ts + \gamma s^2)^{d/2}$, 
the restriction
$\phi_| : \Pp^1 \setminus \{ (\tau_1:1), (\tau_2:1) \} \longrightarrow \bar \CC \setminus \{ P_1,P_2 \}$
is an homeomorphism from $\Pp^1 \setminus \{ (\tau_1:1), (\tau_2:1) \}$ to
 $\CC = \bar \CC \setminus \{ P_1,P_2 \}$. Then $\CC$ is homeomorphic to $\Cc^*$.
\end{proof}

\subsection{Case of $\CC$ being homeomorphic to $\Cc$}

We recall the Abhyankar-Moh-Suzuki \cite{Ar,AM,Ru} theorem:
\begin{theorem} \ 
\label{th:AM}
\begin{enumerate}
 \item \label{it:AM1} Let $t \mapsto (p(t),q(t))$ be an injective polynomial map over $\Cc$
such that the tangent vector $(p'(t),q'(t))$ is never $(0,0)$
then $\deg p$ divides $\deg q$, or $\deg q$ divides $\deg p$.
 \item \label{it:AM2} Let $\CC = (P(x,y)=0) \subset \Cc^2$ be an algebraic plane curve non-singular and homeomorphic 
to $\Cc$ then there exists an algebraic automorphism $\Phi \in \Aut \Cc^2$ such that
$$P \circ \Phi(x,y) = x.$$ 
\end{enumerate}
\end{theorem}

The second statement is the usual form of the Abhyankar-Moh-Suzuki theorem,
it is in fact a consequence of the first (see the proof below), 
for which a more general statement exists \cite{AM}.

\begin{lemma}
\label{lem:AM}
If $\CC = (P(x,y)=k)$, $k$ a generic value, is homeomorphic to $\Cc$ then there exists
$\Phi \in \Aut \Qq^2$ whose inverse has integral coefficients such that
$P \circ \Phi(x,y) = ax+b$, $a,b \in \Qq$.
\end{lemma}

\begin{proof}
The existence of $\Phi \in \Aut \Cc^2$ is Theorem \ref{th:AM}-(\ref{it:AM2}). 
But here we ask the coefficients of $\Phi$ to be rationals
and those of $\Phi^{-1}$ to be integers: we will apply Theorem \ref{th:AM}-(\ref{it:AM1}).
As $\CC$ is homeomorphic to $\Cc$ let $(p(t),q(t))$ be a parametrisation of $\CC$
with $p(t),q(t) \in \Qq[t]$ (it comes from setting $s=1$ in a parametrisation
$\phi = (\bar p(t,s), \bar q(t,s), a s^d)$ of $\bar \CC$).
Then  $\deg p$ divides 
$\deg q$ or $\deg q$ divides $\deg p$.
Suppose that $\delta = \deg p>0$ divides $\deg q$ and
write $p(t) = a_\delta t^\delta + a_{\delta-1} t^{\delta-1} + \cdots$
and $q(t) = b_{\ell \delta} t^{\ell \delta} + b_{\ell \delta-1} t^{\ell \delta -1}+ \cdots$
with $a_i, b_i \in \Qq$ and $\ell \ge 1$.
Write $a_\delta = \frac \alpha \beta$ and $b_{\ell \delta} = \frac{\alpha'}{\beta'}$.
Set the algebraic automorphism of $\Aut \Qq^2$,
$$\Phi_1(x,y) = \bigg(\frac{1}{\alpha \beta} x, \frac{1}{\alpha} \ell y - \frac{\beta^\ell}{\alpha^\ell \beta'} x^\ell\bigg),$$
whose inverse is
$$ \Phi_1^{-1}(x,y) = \bigg(\alpha \beta' x, \alpha' y + \alpha'{\beta'}^{\ell-1}\beta^\ell x^\ell\bigg),$$
whose coefficients are integers.

The composition with $\Phi_1$ yields a parametrisation of $(P-k) \circ \Phi_1(x,y)$ given 
by $p(t), q'(t)$ with $q'(t) \in \Qq[t]$ and $\deg q' < \deg q$.
We repeat this process until one of $p(t)$ or $q(t)$ is a constant
the other being of degree $1$ (because $\CC$ does not have singular points).
Then by the algebraic automorphism $\Phi = \Phi_1 \circ \Phi_2 \circ \cdots$,
whose inverse has integral coefficients, we get
$(P-k) \circ \Phi (x,y) = ax+b$, $a,b \in \Qq$.
\end{proof}

\subsection{Case of $\CC$ being homeomorphic to $\Cc^*$}

We will need the classification over $\Cc$ of polynomials with a generic fiber homeomorphic
to $\Cc^*$, due to W.~Neumann \cite[\S 8]{Ne}.
\begin{theorem}
\label{th:neumann}
If $\CC = (P(x,y)=k)$, $k\neq 0$ a generic value, is homeomorphic to $\Cc^*$ then
$P$ there exists an algebraic automorphism $\Phi \in \Aut \Cc^2$ such that
$$P \circ \Phi(x,y) = x^py^q + \beta$$
$$\text{ or } \quad P \circ \Phi(x,y) =  x^p(yx^r+a_{r-1}x^{r-1} + \cdots + a_0)^q +\beta, $$
with $\beta \in \Cc$, $p >0$, $q >0$, $\gcd(p,q)=1$, $r>0$, $a_0, \ldots, a_{r-1} \in \Cc$ and $a_0 \neq 0$. 
\end{theorem}

We will prove that only some special polynomials of the first type 
can have an infinite number of integral points.

The main result of this part is the following lemma.
\begin{lemma}
\label{lem:neumann}
If $\CC = (P(x,y)=k)$, $k$ a generic value, is homeomorphic to $\Cc^*$ with an infinite number of integral points 
then there exists $\Phi \in \Aut \Qq^2$ such that $P \circ \Phi(x,y) = \alpha(x^2-dy^2)+\beta$, $\alpha \in \Qq^*$, $\beta\in\Qq$.
\end{lemma}

\begin{proof}
By Theorem \ref{th:neumann} we know that the polynomial $P$ has exactly two absolute irreducible factors.
For simplicity of the redaction we suppose in the following that $\beta = 0$. 

\textbf{First case : $P$ is reducible in $\Qq[x,y]$.}

Once again we will prove that $\Phi$ of Theorem \ref{th:neumann} can be chosen with rational coefficients and
its inverse with integral coefficients.
We write $P = \alpha A^pB^q$ the decomposition into irreducible factors with $A,B \in \Qq[x,y]$.
Again for simplicity we suppose $\alpha = 1$.
We will decompose the proof according to the cases of Theorem \ref{th:neumann}. In both cases we see that the curve $(P=0)$ has a non-singular irreducible component homeomorphic to $\Cc$ (the one send by $\Phi^{-1}$ to $(x=0)$). Notice that this curve $(P=0)$
is \emph{not} the curve $\CC$. This component homeomorphic to $\Cc$ is either $(A=0)$ or $(B=0)$, say it is $(A=0)$.
Then, as $A\in\Qq[x,y]$, by the version of Abhyankar-Moh-Suzuki theorem used 
as in Lemma \ref{lem:AM} above, there exists $\Psi \in \Aut \Qq^2$, whose inverse has integral coefficients, such that:
$A \circ \Psi (x,y) = ax+b$, this implies :
$$P \circ \Psi (x,y) = (ax+b)^p Q(x,y)^q.$$

\textbf{Sub-case $P \circ \Phi (x,y) = x^py^q$.}

Then $(Q(x,y)=0)$ is non-singular, homeomorphic to $\Cc$
and the intersection multiplicity with $ax+b$ is $1$.
Then if $(p(t),q(t))$ is a polynomial parametrisation of $(Q(x,y)=0)$
we have $\deg p = 1$. Then as in the proof of Lemma \ref{lem:AM}
by algebraic automorphisms of type $(x,y) \mapsto (\alpha x,\beta y - \gamma x^\ell)$ 
whose  inverse have integral coefficients, we can suppose that $q(t)$ is a constant.
Notice that such automorphisms preserve vertical lines.

Then $Q(x,y)$ becomes $cy+d$, $c,d \in \Qq$, while $ax+b$ remains unchanged.
Then we have find $\Psi' \in \Aut \Qq^2$ whose inverse has integral coefficients such that:
$$P \circ \Psi' (x,y) = (ax+b)^p (cy+d)^q.$$

\textbf{Sub-case $P \circ \Phi (x,y) = x^p(yx^r+a_{r-1}x^{r-1} + \cdots + a_0)^q$.}

$P \circ \Phi (x,y) = x^p(yx^r+a_{r-1}x^{r-1} + \cdots + a_0)^q$ is algebraically equivalent to
$P \circ \Psi (x,y) = (ax+b)^p Q(x,y)^q$  by the algebraic automorphism
$\Phi \circ \Psi^{-1}$. Moreover $\Phi \circ \Psi^{-1}$ should send $x$ to $ax+b$.
Then $\Phi \circ \Psi^{-1}$ is the composition of algebraic automorphisms of type
$(x,y) \mapsto (ax+b,y)$ and $(x,y) \mapsto (\alpha x,\beta y - \gamma x^\ell)$.
This implies that the degree in the variable $y$ remains unchanged.
Then $\deg_y Q(x,y) = \deg_y (yx^r+a_{r-1}x^{r-1} + \cdots + a_0) = 1$.
Then $Q(x,y) = q_1(x)y+q_2(x)$. Due to the asymptotic branches we have $q_1(x,y) = (ax+b)^r$.
And by algebraic automorphisms whose inverse have integral coefficients of type  
$(x,y) \mapsto (\alpha x,\beta y - \gamma x^\ell)$ we can suppose $\deg q_2 < r$.
Then we have found $\Psi' \in \Aut \Qq^2$ with an inverse having integral coefficients such that:
$$P \circ \Psi' (x,y) = (ax+b)^p (y(ax+b)^r+b_{r-1}x^{r-1}+\cdots +b_0)^q,$$
$b_0,\ldots,b_{r-1} \in \Qq$, $b_0 \neq 0$.

\textbf{Conclusion for both sub-cases.}

Now the curve $(P \circ \Psi'(x,y) = k)$ has a finite number of integral points 
since the branches at infinity are asymptotic to horizontal or vertical lines
(with equation $(ax+b=0)$, $(cy+d=0)$ in the first case and $(ax+b=0)$, $(y=0)$ 
in the second case).
Now as $\Psi'^{-1}$ has integral coefficients, an integral point 
$(m,n) \in  (P(x,y)=k) \cap \Zz^2$ is sent to an integral point 
${\Psi'}^{-1}(m,n) \in  (P \circ \Psi'(x,y) = k) \cap \Zz^2$ it implies that
$\CC = (P(x,y)=k)$ also have a finite number of integral points.

\bigskip

\textbf{Second case : $P$ is irreducible in $\Qq[x,y]$.}

Then by Lemma \ref{lem:debes} below it implies that there exist $D,E \in \Qq[x,y]$, $d\in \Zz$
such that $P = C^2-dD^2$. Then $P = (C-\sqrt{d}D)(C+\sqrt{d}D)$ is the decomposition into irreducible factors.

\textbf{Sub-case $P \circ \Phi (x,y) = x^py^q$.}

The by Lemma \ref{lem:debes} we know that $p=1$, $q=1$. And equivalently there exists 
$\Phi' \in \Aut \Cc^2$ such that $P \circ \Phi' (x,y) = (x-\sqrt{d}y)(x+\sqrt{d}y)$.
Then $P \circ \Phi' (x,y) = (C^2-dD^2) \circ \Phi'(x,y) = (x-\sqrt{d}y)(x+\sqrt{d}y)$.
We may suppose that $(C-\sqrt{d}D) \circ \Phi'(x,y) = (x-\sqrt{d}y)$ and
$(C+\sqrt{d}D) \circ \Phi'(x,y) = (x+\sqrt{d}y)$, by addition and subtraction
we get $C \circ \Phi'(x,y) = x$ and $D \circ \Phi'(x,y) = y$.
Then $(CD) \circ \Phi'(x,y) = xy$, with $C,D \in \Qq[x,y]$. As in the first case above
we are now enable to find $\Psi \in \Aut \Qq^2$ (whose inverse as integral coefficients) such that $C \circ \Psi (x,y)=x$,
$D \circ \Psi (x,y)=y$ and $CD \circ \Psi (x,y)=xy$. Now 
$P \circ \Psi(x,y) = (C^2-dD^2) \circ \Psi(x,y) = x^2-dy^2$.

\textbf{Sub-case $P \circ \Phi (x,y) = x^p(yx^r+a_{r-1}x^{r-1} + \cdots + a_0)^q$.}

Again $p=1$, $q=1$, and we may suppose that $(C-\sqrt{d}D) \circ \Phi(x,y) = x$.
We denote $Q = y \circ \Phi^{-1}(x,y)$ i.e.\ $Q \circ \Phi(x,y) = y$.
Then
\begin{align*}
 (C-\sqrt{d}D)(C+\sqrt{d}D) 
     &= P \\
     &= x(yx^r+a_{r-1}x^{r-1} + \cdots + a_0) \circ \Phi^{-1}(x,y) \\
     &= (C-\sqrt{d}D)\big( Q(C-\sqrt{d}D)^r + \cdots \big) \\
\end{align*}
Then 
$$(C+\sqrt{d}D) = \big( Q(C-\sqrt{d}D)^r + \cdots \big).$$
But as $C,D \in \Qq[x,y]$ we have $d = \deg (C+\sqrt{d}D) = \deg (C-\sqrt{d}D)$ and
we get $d = \deg Q + rd$. As $r\ge 1$ we get $\deg Q=0$ which is in contradiction 
with the definition of $Q$. Then this sub-case does not occur.
\end{proof}

\begin{lemma}
\label{lem:debes}
Let $P\in \Qq[x,y]$ such that $P = \alpha A^pB^q$, with $\gcd(p,q)=1$ and with $\alpha \in \Qq^*$, 
$A,B \in\overline{\Qq}[x,y]$ normalized and irreducible (that is to say
$P$ admits exactly two absolute irreducible factors).
Then either $A,B \in \Qq[x,y]$ 
or $p=1$, $q=1$ and there exist $C,D \in \Qq[x,y]$, $d\in \Zz$ non-square such that $P = \alpha(C^2-dD^2)$.
\end{lemma}

The following proof is due to Pierre D\`ebes.
\begin{proof}
Let $a_{i,j} \in \overline{\Qq}$ be the coefficients of $A$.
Let $n$ be the degree of the finite extension $\Qq((a_{i,j}))/\Qq$.
Then there exist exactly $n$ distinct conjugates of $A$.
But for all $\sigma \in \mathrm{Gal}(\overline{\Qq}/\Qq)$, 
$\sigma(A) \in \{A,B\}$. Then $A$ has at most two distinct algebraic conjugates.
Thus $n=1$ or $n=2$.
If $A\notin \Qq[x,y]$ then there exists $a_{i_0,j_0} \notin \Qq$ and 
a $\sigma_0 \in  \mathrm{Gal}(\overline{\Qq}/\Qq)$ such that $\sigma_0(A)=B$.
Then $n=2$ so that the extension  $\Qq((a_{i,j}))/\Qq$ is quadratic. Moreover $p=q=1$. 
This implies the existence of a non-square integer $d$ such that $A,B \in \Qq(\sqrt{d})[x,y]$. 
Now if we write $A = C + \sqrt{d}D$, $C,D \in \Qq[x,y]$ then its algebraic conjugate is $B = C - \sqrt{d}D$. 
\end{proof}

Lemma \ref{lem:AM} and Lemma \ref{lem:neumann} imply Theorem \ref{th:coord} of the introduction.

%%%%%%%%%%%%%%%%%%%%%%%%%%%%%%%%%%%%%%%%%%%%%%%%%%%%%%%%%%%%%%%%%
\section{Number of integral points on polynomials curves}
\label{sec:poly}

\begin{lemma}
\label{lem:supB}
Let $p(t) = a_d t^d + a_{d-1}t^{d-1}+\cdots + a_0 \in \Qq[t]$, $a_d > 0$.
Let $\sigma = - \frac{a_{d-1}}{d a_d}$.
For all $\epsilon > 0$ there exists $B_0 > 0$ such that for all
$B \ge B_0$ and 
$$t_+ = \left( \frac{B}{a_d} \right)^{\frac 1 d} + \sigma + \epsilon,
\quad t_- = -\left( \frac{B}{a_d} \right)^{\frac 1 d} + \sigma - \epsilon,$$
then
$$\forall t \ge t_+ \quad |p(t)| > B \quad \text{ and } \quad \forall t \le t_- \quad |p(t)| > B.$$
\end{lemma}

A similar result holds if $a_d < 0$.

\begin{proof}
We write $t = s+\sigma +\epsilon$ and we look at the asymptotic behaviour
for $p(t)$ when $t$ (and $s$) is large.
\begin{align*}
p(t) &= p(s+\sigma+\epsilon) \\
  &= a_d(s+\sigma+\epsilon)^d + a_{d-1}(s+\sigma+\epsilon)^{d-1}+ \cdots \\
  &= a_d s^d + (d a_d (\sigma+\epsilon) + a_{d-1})s^{d-1} + o(s^{d-1}) \\
  &= a_d s^d + d a_d \epsilon s^{d-1} + o(s^{d-1}). \\
\end{align*}

For $s = \left( \frac{B}{a_d} \right)^{\frac 1 d}$ then $a_d {s}^d = B$ we have
\begin{align*}
p(t_+)
  &= p(s + \sigma+\epsilon) \\
  &= B \cdot \left( 1+ \epsilon d \frac 1 s + o\left( \frac 1 s \right) \right). 
\end{align*}
Then for all sufficiently large $B$ (such that $s>0$ is large enough)
we have $p(t_+) \ge B\left( 1+ \frac{\epsilon d}{2} \frac 1 s \right) $
then $p(t_+) > B$. 
Now the function $t\mapsto p(t)$ is an increasing function for sufficiently large $t$.
Then for all sufficiently large $B$: if $t\ge t_+$ then $p(t) \ge p(t_+) > B.$

Now 
\begin{align*}
p(t_-)
  &= p(- s + \sigma - \epsilon) \\
  &= (-1)^d B \cdot \left( 1 + \epsilon d \frac 1 s + o\left( \frac 1 s \right) \right).
\end{align*}
Then for all sufficiently large $B$, $|p(t_-)| > B$.
And again if $t < t_-$ then $|p(t)| \ge |p(t_-)| > B$.
\end{proof}

For a polynomial $p(t) \in \Qq[t]$ in one variable we defined:
$$M(p,B) = \left\lbrace t \in \Qq \mid p(t) \in \Zz \text{ and } |p(t)| \le B \right\rbrace.$$

\begin{lemma}
\label{lem:supM}
Let $p(t) = \frac 1 b(a_dt^d+\cdots+a_0) \in \Qq[t]$, $a_0,\ldots,a_d,b \in \Zz$, $\gcd(a_0,\ldots,a_d,b)=1$
and $a_d > 0$.
Then for all $\epsilon >0$ there exists $B_0 >0$ such that
for all $B \ge B_0$ we have
$$M(p,B) \le 2a_d^{1-\frac 1 d} b^{\frac 1 d} B^{\frac 1 d} + 1 + \epsilon.$$
\end{lemma}

\begin{proof}
If $t = \frac \alpha \beta \in \Qq$ and $p(\frac \alpha \beta) = k \in \Zz$
then it is well-known that $\beta$ divides $a_d$. Then such $t$ belongs to $\frac 1 {a_d} \Zz$.
Let $\epsilon > 0$ and let $B_0$ be as in Lemma \ref{lem:supB}.
Again by Lemma \ref{lem:supB} if $t>0$ and $|p(t)| \le B$ then 
$t < t_+ = \left( \frac{B}{a_d/b} \right)^{\frac 1 d} + \sigma + \epsilon$.
If $t<0$ and $|p(t)| \le B$ then $|t| < |t_-| = - t_- = \left( \frac{B}{a_d/b} \right)^{\frac 1 d} - \sigma + \epsilon.$
Now the cardinal of $\frac 1 {a_d} \Zz \cap [t_-,t_+]$
is less than $a_d |t_+| + a_d |t_-| + 1 = 
2a_d \left( \frac{B}{a_d/b} \right)^{\frac 1 d} +\sigma - \sigma + 2\epsilon + 1
= 2a_d \left( \frac{bB}{a_d} \right)^{\frac 1 d} + 1 + 2\epsilon$.
\end{proof}

Of course if $p(t)$ is a monic polynomial with integral coefficients, i.e. $b=1$, $a_d=1$, then $M(p,B) \le 2 B^{\frac 1 d} + 1 + \epsilon$.
For example if $p(t)  = t^d$ then $M(p,B) = 2B^{\frac 1 d} + 1$.
The following examples show that the ``$\epsilon$'' is necessary and that the bound of Lemma \ref{lem:supM} is the best one (at least for $a_d=1$).

\begin{example}
\label{ex:epsilon}
Let $p(t) = t^d - 1$ where $d$ is an even number.
Then the following assertion:
$$\exists B_0 >0 \quad \forall B \ge B_0 \quad M(p,B) \le 2B^{\frac 1 d} + 1 \quad \text{is false.}$$
In fact for $k$ any positive integer, set $B_k = p(k) = k^d - 1$.
Then as $d$ is even for all $t \in [-k,k]$ we have $t^d-1 \le k^d-1 = B_k$ then
 $M(p,B_k) = 2k+1$.
Now if the assertion were true we would have
$M(p,B_k) \le 2B_k^{\frac 1 d} + 1$
that it to say $2k+1 \le 2p(k)^{\frac 1 d} + 1$,
then $k^d \le p(k) = k^d-1$ which gives the contradiction.
\end{example}

Here is another example.
\begin{example}
\label{ex:intval}
Let $p(t) = \frac{1}{d!} (t-1)(t-2)\cdots(t-d)$ then
for all $t\in \Zz$ we have $p(t) \in \Zz$. Conversely if
$p(\frac \alpha \beta) \in \Zz$ then $\beta = 1$.
Then for a positive integer $k$ and $B_k = p(k)$, $|p(-k+d)| \le B_k$.
we have $M(p,B_k) \ge 2k+1-d$.
Lemma \ref{lem:supM} applied with $a_d=1$ and $b = d!$
gives $M(p,B) \le 2 (d!)^{\frac 1 d} B_k^{\frac 1 d} + 1 + \epsilon$.
Now $B_k = p(k) \le \frac{k^d}{d!}$
and we get
$$ 2k+1-d  \le M(p,B_k) \le 2k+1+\epsilon.$$
\end{example}

We apply these computations to the situation of our curves.

Let $P(x,y) \in \Qq[x,y]$ be irreducible, let $\CC = (P(x,y) =0)$. Then
$\CC$ is a \defi{polynomial curve} if it admits a polynomial parametrisation
$(p(t),q(t))$, $p(t),q(t) \in \Qq[t]$. Equivalently $\CC$ is a rational 
curve with one place at infinity. Moreover $\deg P = \max(\deg p, \deg q)$.
We will suppose $\deg P = \deg p$ and we write
$p(t) = \frac 1 b(a_dt^d + \cdots + a_0)$ as before.

\begin{lemma}
\label{lem:supN}
Let  $\CC$ be a polynomial curve. Suppose $\deg P = d = \deg p$,
$p(t) = \frac 1 b(a_dt^d + \cdots + a_0)$. Then
for all $\epsilon>0$ there exists $B_0 > 0$ such that for all $B \ge B_0$:
$$N(\CC,B) \le 2a_d^{1-\frac 1 d}b^{\frac 1 d}B^{\frac 1 d} + \frac{(d-1)(d-2)}{2} + 1 + \epsilon.$$
\end{lemma}

The term $\frac{(d-1)(d-2)}{2}$ comes from the number of singular points; for non-singular curves we get the bound of Theorem \ref{th:numb}.

\begin{proof}
An algebraic curve of degree $d$ must have less than $\frac{(d-1)(d-2)}{2}$ singular points, see \cite[p.117]{Fu}.
The other integral points $(p(t),q(t))$ of $\CC$ correspond to rational parameters $t$, see Lemma \ref{lem:param}.
Now we apply Lemma \ref {lem:supM}.
\end{proof}

Examples \ref{ex:epsilon} and \ref{ex:intval}
give polynomials parametrised by $(p(t),t)$
that proves that the bound of Lemma \ref{lem:supN}
is asymptotically sharp.

\enlargethispage{1cm}

%%%%%%%%%%%%%%%%%%%%%%%%%%%%%%%%%%%%%%%%%%%%%%%%%%%%%%%%%%%%%%%%%


\begin{thebibliography}{Mm}
  
\bibitem{Ar} E. Artal-Bartolo, 
\emph{Une d\'emonstration g\'eom\'etrique du th\'eor\`eme d'Abhyankar-Moh.}
 J. Reine Angew. Math.  464  (1995), 97--108.

\bibitem{AM} S.S. Abhyankar, T.T. Moh,  
\emph{Embeddings of the line in the plane.}
J. Reine Angew. Math. 276 (1975), 148--166.

\bibitem{Ay} M. Ayad, 
\emph{Sur les polyn\^omes $f(X,Y)$ tels que $K[f]$ est int\'egralement ferm\'e dans $K[X,Y]$.}
Acta Arith.  105  (2002), 9--28.

\bibitem{Fu} W. Fulton, 
\emph{Algebraic curves.}
Addison-Wesley, reprint of 1969.

\bibitem{Hb} D.R. Heath-Brown, 
\emph{The density of rational points on curves and surfaces.}
 Ann. of Math.  155  (2002), 553--595.

\bibitem{Ma} E. Maillet, 
\emph{D\'etermination des points entiers des courbes alg\'ebriques unicursales \`a coefficients entiers.} 
J. Ecole Polytech. 2, 20 (1919), 115--156.

\bibitem{Ne} W.D. Neumann, 
\emph{Complex algebraic plane curves via their links at infinity.} 
Invent. Math.  98  (1989), 445--489.

\bibitem{Nc} Nguyen Van Chau, 
\emph{Integer points on a curve and the plane Jacobian problem.} 
 Ann. Polon. Math.  88  (2006), 53--58.

\bibitem{PV} D. Poulakis, E. Voskos, 
\emph{On the practical solution of genus zero Diophantine equations.} 
 J. Symbolic Comput.  30  (2000), 573--582. 

\bibitem{Ru} L. Rudolph, 
\emph{Embeddings of the line in the plane.}
J. Reine Angew. Math. 337 (1982), 113--118. 

\bibitem{Se} J.-P. Serre, 
\emph{Lectures on the Mordell-Weil theorem.}   
Aspects of Mathematics, Vieweg 1997.

\bibitem{Wa} Y. Walkowiak, 
\emph{Th\'eor\`eme d'irr\'eductibilit\'e de Hilbert effectif.}
Acta Arith.  116  (2005), 343--362.

\end{thebibliography}
\end{document}